\documentclass[a4paper,12pt]{article}

\usepackage{amsmath,amsfonts,latexsym,theorem}
\nonstopmode

\numberwithin{equation}{section}

\newtheorem{Theorem}{Theorem}[section]
\newtheorem{Definition}{Definition}[section]
\newtheorem{Proposition}{Proposition}[section]
\newtheorem{Lemma}{Lemma}[section]

\newtheorem{Example}{Example}[section]
\newtheorem{Corollary}{Corollary}[section]
\newenvironment{Proofc}[1]{\smallskip\par\noindent\textsc{#1}\quad}%
  {\hfill$\Box$\bigskip\par}
\newenvironment{Proof}{\begin{Proofc}{Proof}}{\end{Proofc}}

\newtheorem{Remark}{Remark}[section]

\def\a{\alpha}
\def\b{\beta}
\def\d{\delta}

\def\g{\gamma}
\def\G{\Gamma}
\def\l{\lambda}

\def\o{\omega}
\def\th{\theta}

\def\e{\varepsilon}

\def\pd{\partial}

\newcommand{\hx}{{\hat x}}

\newcommand{\p}{\partial}

\newcommand{\R}{{\mathbb R}}
\newcommand{\N}{{\mathbb N}}

\def\pd{\partial}

\begin{document}
\title{The vanishing viscosity limit for Hamilton-Jacobi equations on  networks}

\author{Fabio Camilli\footnotemark[1] \and Claudio Marchi\footnotemark[2] \and Dirk Schieborn\footnotemark[3]}

\date{version: \today}
\maketitle

\footnotetext[1]{Dip. di Scienze di Base e Applicate per l'Ingegneria,  ``Sapienza'' Universit{\`a}  di Roma, via Scarpa 16,
 00161 Roma, Italy, ({\tt e-mail:camilli@dmmm.uniroma1.it})}
\footnotetext[2]{Dip. di Matematica, Universit\`a di Padova, via Trieste 63, 35121 Padova, Italy ({\tt marchi@math.unipd.it}).}
\footnotetext[3]{Eberhard-Karls University,  T\"ubingen, Germany  ({\tt e-mail:Dirk@schieborn.de})}

\begin{abstract}
For a Hamilton-Jacobi equation defined on   a network, we introduce its vanishing viscosity approximation. The elliptic equation is given on the  edges and  coupled with Kirchhoff-type conditions at the transition vertices.
We prove that there exists exactly one solution of this elliptic approximation and mainly that, as the viscosity vanishes, it converges to the unique solution of the original problem.

\end{abstract}
 \begin{description}
\item [{\bf MSC 2000}:] 35R02, 35F21,35D40, 35B40.
\item [{\bf Keywords}:] Vanishing viscosity, network, Hamilton-Jacobi equation, viscosity solution, maximum principle.
\end{description}
%
%
%
\section{Introduction}\label{intro}
The study of partial differential equations on  networks  arise is several applications as information networks (internet, social networks, email exchange), economical networks (business relation between companies, postal delivery and traffic routes),  biological networks (neural networks, food web, blood vessel, disease transmission).

Starting with the seminal work of Lumer \cite{lu}, a fairly complete theory for linear and semilinear equations on networks
has been developed in the last 30 years (for instance, see: Lagnese et al. \cite{lls}, Von Below et al. \cite{vbn}, Engel et al. \cite{ekns}, Freidlin et al. \cite{fsh,fw}). Only in   recent times it has been initiated the study of some classes of fully nonlinear equations,  such as  conservation laws (see \cite{cg,gp} and reference therein) or   Hamilton-Jacobi equations (see \cite{acct,csm, ghn,imz,sc}).

All the     approaches to Hamilton-Jacobi equations aim to extend the concept of viscosity solution (see \cite{b, bcd}) to networks, but they differ  for the assumptions made on the Hamiltonians at the vertices. Hence, different frameworks reflect in different definitions of viscosity solutions, even if all of them  give  existence and uniqueness of the solution.
However, any generalization of viscosity solution should preserve the other main  features of existing theory such as  stability with respect to uniform convergence and the method of vanishing viscosity.

In this paper we aim to show that the definition of solution introduced in \cite{sc} is consistent with vanishing viscosity method, which
consists in approximating the original nonlinear problem by a family of semilinear ones.
The difficulty is thus transferred to the question, whether the approximating family of solutions converges.

The first step establishes existence and uniqueness of classical solutions to the viscous Hamilton-Jacobi equation on networks. In doing so, the necessity of an extra condition at transition vertices   becomes clear.  We impose   the classical Kirchhoff condition which establishes  a relation among the outer normal derivatives of the solution along the   edges incident the same vertex. The Kirchhoff condition can be thought of as an extension of the  ``averaging effect'' of the viscosity term on the vertices.

The second step is to prove some a priori estimates, uniform in the viscosity parameter. These estimates are obtained by explicit arguments which take advantage of the intrinsic one dimensional nature of the problem.

The final step is the convergence of the solution of viscous approximation to the one of the starting problem. Obviously this issue requires a special care at the vertices, while it follows by classical arguments inside the edges.

The paper is organized as follows.
In Section \ref{s1}, we introduce some notations, the standing assumptions and recall the definition of viscosity solution.
In Section~\ref{s2} we  study existence and uniqueness of the solution to  the second order problem. Section~\ref{s3}  is devoted to the proof of the  a priori estimates, whereas in Section~\ref{s4} we show the convergence of the vanishing viscosity method; we work out in detail the eikonal problem in Section~\ref{esempio}. In Appendix~\ref{AppendixB} we prove some technical lemmas.


%
\section{Notations and preliminary definitions}\label{s1}
\subsection{Topological network}
A topological network   is a collection of points in $\R^n$ connected by continuous, non self-intersecting curves.
More precisely (see \cite{lu, sc}):
\begin{Definition}\label{s1:D1}let $V=\{v_i,\,i\in I\}$ be a finite collection of points in $\R^n$ and
let $\{\pi_j,\,j\in J\}$ be a finite collection of smooth, non self-intersecting curves in $\R^n$ given by
$\pi_j:[0,l_j]\to\R^n,\, l_j>0$. For $e_j:=\pi_j((0,l_j))$ and $\bar e_j:=\pi_j([0,l_j])$, 
assume that
\begin{itemize}
  \item[i)] $\pi_j(0), \pi_j(l_j)\in V$, and
$\#(\bar e_j\cap V)=2$ for all $j\in J$,
  \item[ii)] $\bar e_j\cap \bar e_k\subset  V$, and $\#(\bar e_j\cap \bar e_k)\le 1$ for all $j,k \in J$, $j\neq k$.
  \item[iii)] For all $v, w \in V$ there is a path with end-points  $v $ and $w$ (i.e.  a sequence of edges $\{e_j\}_{j=1}^N$ such that
  $\#(\bar e_j\cap \bar e_{j+1})=1$ and  $v\in \bar e_1$, $w\in \bar e_N$).
\end{itemize}
Then $\G:=\bigcup_{j\in J}\bar e_j\subset \R^n$ is called a (finite)  \emph{topological network} in $\R^n$.
\end{Definition}
In the following we always identify   $x\in \bar e_j$ with   $y=\pi_j^{-1}(x)\in [0,l_j]$.
For $i\in I$ we set $Inc_i:=\{j\in J:\,e_j \,\text{is incident to}\,v_i\}$, moreover two vertices $v_i$, $v_j$ are said adjacent (in symbols
$v_i\, \mathrm{adj} \,v_j$) if there exists $k\in J$ such that $v_i,\, v_j\in e_k$.\\
Observe that the parametrization of the arcs $e_j$ induces an orientation which can be expressed by the \emph{signed incidence matrix} $A=\{a_{ij}\}$ with
\begin{equation}\label{1:1}
   a_{ij}:=\left\{
            \begin{array}{rl}
              1 & \hbox{if $v_i\in\bar e_j$ and $\pi_j(0)=v_i$,} \\
              -1 & \hbox{if $v_i\in\bar e_j$ and $\pi_j(l_j)=v_i$,} \\
              0 & \hbox{otherwise.}
            \end{array}
          \right.
\end{equation}
Given a nonempty set $I_B\subset I$, we define  $\partial \G:=\{v_i,\,i\in I_B\}$; we assume $i\in I_B$ whenever $\# Inc_i=1$ (see Remark~\ref{rmk-A} below).
For $I_T:=I\setminus I_B$, we call $\{v_i:\,i\in I_B\}$ the set of boundary vertices and $\{v_i:\,i\in I_T\}$ the set of transition vertices.

\subsection{Function spaces}
For any function $u:  \G\to\R$ and each $j\in J$ we denote by $u^j:[0,l_j]\to \R$ the restriction of $u$ to $\bar e_j$, i.e.  $u^j(y)=u(\pi_j(y))$ for $y\in [0,l_j]$.
For $\a \in\N$, we define differentiation along an edge $e_j$ by
\[\pd^\a_ju(x):= \frac{d^\a u^j}{d y^\a}  (y),\qquad\text{for $y=\pi^{-1}_j(x)$, $x\in e_j$}\]
and at a vertex $v_i$ by
\[\pd^\a_ju(v_i):= \frac{d^\a u^j}{d y^\a}  (y)\qquad\text{for $y=\pi^{-1}_j(v_1)$,  $j\in Inc_i$.}\]
\begin{Definition}\label{def:USC}

\begin{itemize}
\item[i)] We say that a function  $u$ belongs to $USC( \G)$ (respectively, to $LSC( \G)$) if it is upper (resp., lower) semicontinuous with respect to the topology induced by $\R^n$ on $ \G$. In other words, $u\in USC( \G)$ if and only if $u^j\in USC([0,l_j])$ for every $j\in J$ and $u^j(\pi_j^{-1}(v_i))=u^k(\pi_k^{-1}(v_i))$ for every $i\in I$, $j,k\in Inc_i$; an analogous property holds for $u\in LSC(\G)$.
\item[ii)] We say that a function  $u$ is  continuous in $\G$ and we write $u\in C(\G)$  if it is continuous with respect to the subspace topology of $\G$, namely, $u^j\in C([0,l_j])$ for any $j\in J$ and $u^j(\pi_j^{-1}(v_i))=u^k(\pi_k^{-1}(v_i))$  for any $i\in I$, $j,k\in Inc_i$.
\item[iii)] We say that $u\in C^k(\G)$ if $u\in C(\G)$ and if $u^j\in C^k([0,l_j])$ for $j\in J$.
\item[iv)] For any collection $\b=(\b_{ij})_{i\in I_T,\,j\in Inc_i}$ with $\b_{ij}\geq 0$, we say that $u\in C^k_{*,\b}(\G)$ if $u\in C^k(\G)$, $k\ge 1$, and there holds
\begin{equation}\label{Kir}
S^i_\b u:= \sum_{j\in Inc_i}\b_{ij}a_{ij}\pd_j u(v_i)=0 \qquad \forall i\in I_T.
\end{equation}
\end{itemize}
\end{Definition}
\begin{Remark}
Condition \eqref{Kir} is known in the literature as the Kirchhoff condition. In a way,  differentiability
of a function along the edges means that the slopes in outward (or inward)
direction with respect to each given point add up to zero. At vertices, this condition
naturally generalizes to the Kirchhoff condition.
\end{Remark}
\subsection{Viscosity solutions}
A  Hamiltonian $H:\G\times \R\times \R\to \R$ is a collection of operators $(H^j)_{j\in J}$ with $H^j: [0,l_j]\times\R\times\R\to \R$.
Along the paper  we will consider the following  conditions
\begin{align}
    &H^j\in C^0([0,l_j]\times\R\times \R), \quad j\in J;\label{1:H0}\\
    &H^j(x,\cdot,p)\quad \text{is nondecreasing for all $(x,p)\in [0,l_j]\times \R$,}\quad j\in J;\label{1:H1}\\
   &H^j(v_i,r,\cdot)\quad \text{is nondecreasing in $(0,+\infty)$ for any $i\in I_T$, $r\in\R$;}\label{1:H3}\\
   &H^j(x,r,\cdot)\to +\infty\quad \text{as $|p|\to \infty$ uniformly in }(x,r)\in [0,l_j]\times[-R,R], j\in J;\label{1:H4}\\
   &H^j(\pi_j^{-1}(v_i),r,p)=H^k(\pi_k^{-1}(v_i),r,p)\,\text{for any  $r\in \R$, $p\in \R$, $i\in I_T$, $j,k \in Inc_i$;}\label{1:H5}\\
   &H^j(\pi_j^{-1}(v_i),r,p)=H^j(\pi_j^{-1}(v_i),r,-p)\,\text{for any $r\in \R$, $p\in \R$, $i\in I_T$, $j\in Inc_i$.}\label{1:H6}
\end{align}
\begin{Remark}\label{1:rem1}
 Assumptions~\eqref{1:H5}-\eqref{1:H6} represent compatibility conditions of $H$  at the vertices of  $\G$, i.e.  continuity at the
vertices and independence of the orientation of the incident arc, respectively (the network is not oriented).
\end{Remark}
\begin{Example}
The operator $H(x,r,p):=|p|^\a +b(x)r+f(x)$ satisfies assumptions~\eqref{1:H0}-\eqref{1:H6} provided that $\a>0$, $b,f\in C^0(\G)$ and $b(x)\geq 0$ for every $x\in\G$.
\end{Example}

On the graph~$\G$, we consider the Hamilton-Jacobi equation
\begin{equation}\label{HJ}
  H(x,u, \pd u)=0,\qquad x\in\G,
\end{equation}
namely, on each edge $e_j$, we address the Hamilton-Jacobi equation
\[
  H^j(y,u^j(y), \pd_j u)=0,\qquad y\in [0,l_j].
\]
In the next definitions we introduce the class of test functions and solution of \eqref{HJ}.
\begin{Definition}\label{1:def2}
Let $\phi\in C(\G)$.
\begin{itemize}
  \item[i)] Let $x\in e_j$, $j\in J$. We say that $\phi$ is test function at $x$, if $\phi^j$ is differentiable at $\pi_j^{-1}(x)$.
  \item[ii)] Let $ x=v_i$, $i\in I_T$, $j,k\in Inc_i$, $j\neq k$. We say that $\phi$ is $(j,k)$-test function at $x$, if
  $\phi^j$ and $\phi^k$ are differentiable at $ \pi_j^{-1}(x)$ and $ \pi_k^{-1}(x)$, respectively  and
\begin{equation}\label{1:2}
a_{ij}\pd_j \phi(\pi_j^{-1}(x)  )+a_{ik}\pd_k \phi(\pi_k^{-1}(x) )=0,
\end{equation}
where $(a_{ij})$ as in \eqref{1:1}.
\end{itemize}
\end{Definition}
\begin{Definition}\label{1:def3}
A function $u\in\text{USC}( \G)$ is called a (viscosity) subsolution of \eqref{HJ} in $\G$ if the following holds:
\begin{itemize}
  \item[i)] If  $x\in e_j$, $j\in J$,    for any test function  $\phi$     for which
  $u-\phi$ attains a local maximum at $x$, we have
 \[H^j(\pi_j^{-1}(x),u^j(\pi_j^{-1}(x)), \pd_j \phi (\pi_j^{-1}(x)))\le 0.\]
\item [ii)] If $x=v_i$, $i\in I_T$,  for any $j,k\in Inc_i$ and any $(j,k)$-test function  $\phi$     for which
  $u-\phi$ attains a local maximum at $x$ relatively to $\bar e_j\cup \bar e_k$, we have
 \[H^j(\pi_j^{-1}(x),u^j(\pi_j^{-1}(x)), \pd_j \phi(\pi_j^{-1}(x)) )\le 0.\]
\end{itemize}
 A function $u\in\text{LSC}(\G)$ is called a (viscosity) supersolution of \eqref{HJ} in $\G$ if the following holds:
\begin{itemize}
  \item[i)] If  $x\in e_j$, $j\in J$,    for any test function  $\phi$     for which
  $u-\phi$ attains a local minimum at $x$, we have
 \[H^j(\pi_j^{-1}(x),u^j(\pi_j^{-1}(x)), \pd_j \phi(\pi_j^{-1}(x)) )\ge 0.\]
  \item [ii)] If $x=v_i$, $i\in I_T$, for any  $j\in Inc_i$, there exists $k\in Inc_i$, $k\neq j$, (said $i$-feasible for $j$ at $x$) such that for any  $(j,k)$-test function  $\phi$
   for which $u-\phi$ attains a local minimum at $x$ relatively to $\bar e_j\cup \bar e_k$, we have
 \[ H^j(\pi_j^{-1}(x),u^j(\pi_j^{-1}(x)), \pd_j \phi(\pi_j^{-1}(x)) )\ge 0.\]
\end{itemize}
A continuous function $u\in C(\G)$ is called a (viscosity) solution of  \eqref{HJ}  if it is both a viscosity subsolution and a viscosity supersolution of \eqref{HJ}.
\end{Definition}
\begin{Remark}
It is important to observe that the definitions of subsolution and supersolution are not symmetric at the vertices. As observed
 in \cite{sc} for the equation $|\pd u|^2=1$, a definition of supersolution similar to the one of subsolution would not characterize
 the correct solution, i.e. the distance from the boundary.
\end{Remark}
\begin{Remark}\label{rmk-A}
The definition of solution does not involve the vertices $v_i\in\p\G$: at these points no  ``transition'' condition is required.
Wlog, we assume $\# Inc_i=1$ for any $i\in I_B$. Actually, whenever $i\in I_B$ and $\# Inc_i>1$, the problem is equivalent to the one obtained by splitting the common endpoints of the edges incident  $v_i$.
\end{Remark}
\subsection{Perron method and comparison principle}

In this section we collect some results on the well posedness of the
Hamilton-Jacobi equations \eqref{HJ}.
Concerning the existence of a solution we have the following result; for the proof, obtained via  Perron's method, we refer the reader to~\cite[Thm6.1]{csm}.
\begin{Theorem}\label{Perron}
Assume \eqref{1:H0}-\eqref{1:H6}  and that there is a viscosity subsolution $w\in USC(\G)$ and a viscosity supersolution $W\in LSC(\G)$ of \eqref{HJ} such that $w\leq W$ and $w_*(x)=W^*(x)=g(x)$ for $x\in\pd \G$.
Let the function $u:\G\to\R$ be defined by $u(x):=\sup_{v\in X} v(x)$ where
\[X=\{v\in USC(  \G):\, \text{$v$ is a viscosity subsolution of \eqref{HJ} with $w\le v\le W$ on $\G$}\}.\]
Then, $u^*$ and $u_*$ are respectively  a sub- and a supersolution to problem~\eqref{HJ} with $u=g$ on $\pd \G$.
\end{Theorem}
The  proof of the following theorem relies  on the classical doubling of variable argument;
for  the detailed proof, we refer to \cite[Thm5.1]{sc} and to \cite[Lem5.2]{s}.
\begin{Theorem}\label{s1:T2}
Assume \eqref{1:H0}-\eqref{1:H6}.
\begin{itemize}
\item[$(a)$]
Assume
\begin{equation}\label{Hr}
\text{$ H^j(y,\cdot,p)$ is strictly increasing for any $y\in [0,l_j]$, $p\in \R$, $j\in J$.}
\end{equation}
Let $u_1$ and $u_2$ be  respectively a bounded super- and a bounded  subsolution of \eqref{HJ} such that $u_1(v_i) \ge u_2(v_i)$ for all $i\in I_B$.
Then $u_1\ge u_2$ in $\G$.
\item[$(b)$] Let $u_1$ and $u_2$ be  respectively a supersolution to \eqref{HJ} and a subsolution to
\[H(x,u,\p u)=g(x)\qquad x\in\G\]
with $g\in C(\G)$, $g<0$. Then there holds $u_1\ge u_2$ in $\G$, provided that $u_1(v_i) \ge u_2(v_i)$ for all $i\in I_B$.
\end{itemize}
\end{Theorem}
Finally, let us state a stability result (see \cite[Prp3.2]{sc}):
 \begin{Proposition}\label{s1:P2}
 Assume \eqref{1:H0}-\eqref{1:H6}.
Let $u_n$ be a solution of
\[
   H_n(x,u_n,\pd u_n)=0,\qquad x\in \G, n\in \N.
\]
Assume that, as $n\to\infty$, $H_n(x,r,p)\to H(x,r,p)$ locally uniformly
and $u_n\to u$ uniformly in $\G$. Then $u$ is  a solution of  \eqref{HJ}.
\end{Proposition}
\begin{Remark}
For the  Hamilton-Jacobi equation \eqref{HJ} it is well known  that a smooth solution
  will not exist in general. Furthermore it is equally easy to see  that the Kirchhoff
condition \eqref{Kir} is not   satisfied. Continuity is the
only property of a  solution to \eqref{HJ} which is reasonable to expect.
\end{Remark}
\section{The viscous eikonal equation on networks}\label{s2}
In this section we study the existence and the uniqueness of a classical solution  to second order equations coupled with Kirchhoff condition. 
\subsection{Linear problems}
We consider the following  class of linear problems on $\G$
 \begin{equation}\label{s2:linear}
L^jw(x)+g^j(x)=0\quad\text{$x\in e_j$,\,$j\in J$,}\qquad
w(v_i)=\g_i\quad \in I_B.
\end{equation}
where $L=(L^j)_{j\in J}$ is a collection of elliptic linear operators of the form
\begin{align}
      &L^j w(x):=a^j(x)\pd_j^2 w(x) +b^j(x)\pd_j w(x)-c^j(x) w(x)\qquad   x\in e_j,\,j\in J.\label{s2:1}
\end{align}
We assume the following hypotheses
\begin{equation}\label{2:H0}
a^j, b^j, c^j, g^j\in C([0,l_j]),\quad a^j(x)\ge \l>0,\quad c^j(x)\ge 0\quad\forall x\in [0,l_j],\, j\in J
\end{equation}

Let us now state a maximum principle for problem~\eqref{s2:linear}.
\begin{Theorem}\label{s2:T1}
Let $L=(L^j)$, $S_\b=(S^i_\b)_{i\in I}$ be defined as in \eqref{s2:1}-\eqref{2:H0}
 and respectively in \eqref{Kir} with $\sum_{j\in Inc_i} \b_{ij}>0$ for each $i\in I_T$.
Assume that the function $w\in C^2(\G)$ satisfies
 \begin{equation}\label{s2:3}
 L^jw(x)\ge 0 \quad\,x\in e_j, j\in J\quad\text{and}\quad S^i_\b w\ge 0\quad i\in I_T.
 \end{equation}
Then $w$ attains a nonnegative maximum in $\G\setminus \p\G$ if, and only if, it is constant.
A similar result  holds for the minimum of $w$ if we revert the inequalities in \eqref{s2:3}.
\end{Theorem}
\begin{Proof}
We set $M:=\max w$ and $A:=\{x\in\G\setminus \pd \G:\,w(x)=M\}$. We proceed by contradiction assuming $M\geq 0$ and $A\ne \emptyset$. For the sake of clarity, we split the arguments in two cases.

{\it Case (I).} We assume that $Lw>0$, $S_\b w>0$
 and $x_0\in A$.
 If $x_0\in e_j$ for some $j\in J$, then we have: $\pd_j w(x_0)=0$ and $\pd_j^2 w(x_0)\le  0$, a contradiction to $L^j w>0$.  If $x_0=v_i$ for some $i\in I_T$, then we have $a_{ij}\pd_j w(v_i) \le 0$ for all $j\in Inc_i$, hence $S^i_\b w \le 0$, a contradiction.


{\it Case (II).} We assume that $Lw\ge 0$, $S_\b w\ge 0$ and $x_0\in A$. By the continuity of $w$, one of the following two cases must occur somewhere in $\G$
\begin{itemize}
\item[i)] for some $j\in J$, $x_0\in e_j$ and $w(y)<w(x_0)$ for some $y\in e_j$,
\item[ii)] for some $j\in J$, $x_0=v_i$ and $w(y)<w(x_0)$ for some $y\in e_j$ with $j\in Inc_i$.
\end{itemize}
In case $(i)$, the (nonconstant) function $w^j$ solves $L^j w^j\geq 0$ in $(0,l_j)$ and it attains a nonnegative maximum inside $(0,l^j)$. This situation is impossible by classical results (see \cite[Ch.1]{PW}).

Let us consider case $(ii)$. Now it suffices to prove the statement in the network~$\G_0:=\cup_{j\in Inc_i} \bar e_j$. Moreover, wlog, we shall assume $\pi_j(0)=v_i$ for any $j\in Inc_i$ and $y\in e_{\bar \j}$.
We claim that there exists a function $\phi\in C^2(\G_0)$ such that
\begin{equation}\label{MP:phi}
L^j\phi^j(x)>0 \quad \forall x\in e_j,\, j\in Inc_i, \qquad S^i_\b \phi>0,\qquad \phi\geq 0, \qquad \phi(v_i)=0.
\end{equation}
To this end, we define $\phi^j(x):= e^{\a_j x}-1$ (for $j\in Inc_i$)
with a parameter $\a_j$ such that $L^j\phi^j>0$. In order to have this inequality, it suffices to choose $\a_j>0$ such that there holds
\[\l \a_j^2-\|b^j\|_\infty \a_j -\|c^j\|_\infty>0.\]
Moreover, we have: $S^i_\b \phi=\sum_{j\in Inc_i}\b_{ij} \p_j \phi^j=\sum_{j\in Inc_i}\b_{ij}\a_j>0$. Hence, our claim \eqref{MP:phi} is completely proved.

Fix $\eta:=(w(v_i)-w(y))(e^{\a_{\bar \j} l_{\bar \j}}-1)^{-1}$ (note $\eta >0$ by our assumptions) and introduce the function $\tilde w(x):=w(x)+\eta \phi(x)$, $x\in \G_0$. We observe that there holds
\begin{align*}
&S^i_\b \tilde w=S^i_\b  w +\eta S^i_\b \phi>0,\qquad L^j \tilde w^j=L^j w^j+\eta L^j\phi^j>0 \quad \forall j\in Inc_i\\
&\tilde w(v_i)=w(v_i), \qquad
\tilde w(y)=w(y)+(w(v_i)-w(y))\frac{\phi(y)}{e^{\a_{\bar \j} l_{\bar \j}}-1}<w(v_i).
\end{align*}
Invoking case $(i)$ we obtain a contradiction.
\end{Proof}


\begin{Theorem}\label{linear:exi}
There exists a unique solution $u\in C^2_{*,\b}(\G)$ to problem \eqref{s2:linear}.
\end{Theorem}
\begin{Proof}
By standard arguments (see~\cite[Ch.1]{PW}), uniqueness is an immediate consequence of Theorem~\ref{s2:T1}.
Existence of a solution to \eqref{s2:linear} is proved in \cite[Thm3.3]{fw} (see also the related comments and \cite{fsh})
via a probabilistic representation formula. In fact a solution of \eqref{s2:linear} can be represented as
\[u(x)={\mathbb E}_x\{\int_0^\tau e^{-c(Y(s))}g(Y(s))ds+ e^{-c(Y(\tau))}\g_{i(\tau)}\}\]
where $Y(s)$ is a Markov process defined on the graph which on each edge $e_j$ solves the stochastic differential equation
\[dY(s)=b^j(Y(s))ds+a^j(Y(s))dW(s),\]
 $\tau=\inf\{t>0:\,Y(t)\in\pd\G\}$ and $i(\tau)\in I_B$ is  such that $Y(\tau)=v_{i(\tau)}\in  \pd \G$. In this interpretation the Kirchhoff condition
 \eqref{Kir} implies that the process almost surely spends zero time at each transition vertex $v_i$, (see \cite[Thm3.1]{fw}) while the term $ \b_{ij}/(\sum_{j\in Inc_i}\b_{ij})$ is the probability that $Y(t)$ enters in the edge~$e_j$ when it is in $v_i$.
\end{Proof}
\subsection{Semi-linear problems}
\begin{Theorem}\label{Teo:logar}
For any $\e>0$, there exists a unique solution  $u_\e\in C^2_{*,\beta}(\G)$ of
\begin{equation}\label{HJe}
-\e \pd_j^2u+|\pd_j u|^2-f(x)=0\quad  x\in e_j,\,j\in J ,
\qquad u(v_i)=g_i,\quad i\in I_B
\end{equation}
where $f$ is a continuous, non negative function on $\G$.
\end{Theorem}
\begin{Proof}
We consider  the logarithmic transformation (see \cite{fs}): $u_\e=-\e\ln(w_\e+1)$.
Invoking Theorem~\ref{linear:exi}, we have that for any $\e>0$ there exists a unique solution $w_\e\in C^2_{*,\beta}(\G)$ to the linear problem
\[
\e^2\pd_j^2 w_\e-f(x)w_\e- f(x)=0\quad\quad  x\in e_j,\,j\in J,
\qquad w_\e(v_i)=e^{-\frac{g_i}{\e}}-1,\quad i\in I_B
\]
Hence, reversing the logarithmic transformation, we conclude that there exists a unique solution to \eqref{HJe}.
\end{Proof}

Another consequence of Theorem \ref{s2:T1} is the following comparison principle
\begin{Corollary}\label{s2:C2}
Assume that $H=(H^j)_{j\in J}$ satisfies \eqref{1:H0}-\eqref{1:H1} and
\begin{equation}\label{s2:reeegular}
\text{$ H^j(x,\cdot,\cdot)\in C^1(\R\times \R)$ for any $x\in (0,l_j)$, $j\in J$.}
\end{equation}
Let $w_1, w_2\in C^2(\G)$ be such that
\begin{equation}\label{s2:CP1}
\left\{\begin{array}{ll}
-\e \pd^2_j w_1+H^j(x,w_1,\pd_jw_1)\ge -\e \pd^2_j w_2+H^j(x,w_2,\pd_jw_2) &x\in e_j,\,j\in J,\\
S^i_\b w_1 \le S^i_\b w_2  &i\in I_T\\
w_1(v_i)\ge w_2(v_i)  &i\in I_B
\end{array}\right.
\end{equation}
Then $w_1\ge w_2$ on $\G$.
\end{Corollary}
\begin{Proof}
Set $A=\{w_2>w_1\}\subset\G$; the function~$w:=w_2-w_1$ is a solution to
\[
\left\{\begin{array}{ll}
 \e\pd_j^2 w+\tilde b^j(x)\pd_j w -\tilde c^j(x)w\ge 0 \qquad   &x\in e_j\cap A,\,j\in J,\\
    S^i_\b w\ge 0&  i\in I_T\cap A\\
    w(v_i)\le 0&i\in I_B\cap A
   \end{array}
   \right.
\]
where
\begin{align*}
\tilde b^j(x)&=-\int_0^1\frac{\partial H^j}{\partial p}\left(x, \th w_1+(1-\th)w_2, \th \pd_j w_1+(1-\th)\pd_j w_2\right)d\th\\
\tilde c^j(x)&=\int_0^1\frac{\partial H^j}{\partial r}\left(x, \th w_1+(1-\th)w_2, \th \pd_j w_1+(1-\th)\pd_j w_2\right)d\th.
\end{align*}
By Theorem \ref{s2:T1}, $w$ cannot attain a local nonnegative maximum inside the open set $A$. As we have $A\cap\pd\G=\emptyset$, it follows that $A$ is empty and $w_1\ge w_2$ in $\G$.
\end{Proof}


\subsection{Other comparison principles for \eqref{s2:CP1}}

For the sake of completeness, we establish some comparison principles for problem~\eqref{s2:CP1} under assumptions different from Corollary \ref{s2:C2}; especially, in both of them we shall drop the regularity condition~\eqref{s2:reeegular}.
In the former we require the strict monotonicity of $H$ with respect to $u$, while in the latter we require a linear growth of $H$ with respect to $u$ and $\pd u$.
\begin{Proposition}\label{s2:C2bis}
Assume that $H=(H^j)_{j\in J}$ satisfies \eqref{1:H0}-\eqref{1:H1} and \eqref{Hr}.
Let the functions $w_1, w_2\in C^2(\G)$ satisfy \eqref{s2:CP1} with $\b_{ij}>0$ for any $i\in I_T$, $j\in Inc_i$. Then $w_1\ge w_2$ on $\G$.
\end{Proposition}

\begin{Proof}
We argue by contradiction assuming $\max_{\G}(w_2-w_1)=:\d>0$.
Let $x_0$ be a point where $w_2-w_1$ attains its maximum; whence $x_0\in \G$. The point~$x_0$ either belongs to some edge or it coincides with a transition vertex.
Assume that, $x_0$ belongs to some edge $e_j$. By their regularity, the functions $w_1$ and $w_2$ fulfill
\[
w_2(x_0)=w_1(x_0)+\d,\qquad \pd_j w_2(x_0)=\pd_j w_1(x_0),\qquad \pd_j^2 w_2(x_0)\leq \pd_j^2 w_1(x_0).\]
In particular, we deduce
\begin{align*}
-\e \pd_j^2w_1(x_0)+H(x_0, w_1(x_0),\pd_j w_1(x_0))&\leq
-\e \pd_j^2w_2(x_0)+H(x_0, w_2(x_0)-\d,\pd_j w_2(x_0))\\ &<
-\e \pd_j^2w_2(x_0)+H(x_0, w_2(x_0),\pd_j w_2(x_0))
\end{align*}
which contradicts the first relation in~\eqref{s2:CP1}.

Assume that $x_0=v_i$ for some $i\in I_T$. Being regular, the functions $w_1$ and $w_2$ fulfill $a_{ij}\pd_j w_2(v_i)\leq a_{ij}\pd_j w_1(v_i)$. We claim $\pd_j w_2(v_i)= \pd_j w_1(v_i)$ for each $j\in Inc_i$. In order to prove this equality we proceed by contradiction and we assume that $a_{ij}\pd_j w_2(v_i)< a_{ij}\pd_j w_1(v_i)$ for some $j\in Inc_i$. In this case we get $S^i_\b w_2<S^i_\b w_1$ which contradicts the second hypothesis in \eqref{s2:CP1}; therefore, our claim is proved.
Moreover, since $w_1(x_0)=w_2(x_0)-\d$, we deduce
\[H(x_0, w_1(x_0),\pd_j w_1(x_0))=H(x_0, w_2(x_0)-\d,\pd_j w_2(x_0))<
H(x_0, w_2(x_0),\pd_j w_2(x_0)).\]
Taking into account the regularity of $H$ and of $w_i$ ($i=1,2$), we infer that in a sufficiently small neighborhood $B_\eta(v_i)$ there holds
\[H(x, w_1(x),\pd w_1(x))<H(x, w_2(x),\pd w_2(x))\]
This inequality and the first relation in \eqref{s2:CP1} entail
\[
\e \pd^2_j(w_2-w_1)\geq H(x, w_2(x),\pd_j w_2(x))-H(x, w_1(x),\pd_j w_1(x))>0
\]
which, together with $\pd_j w_2(v_i)= \pd_j w_1(v_i)$, contradicts that $w_2-w_1$ attains a maximum in $v_i$.
\end{Proof}

\begin{Proposition}\label{s2:C2ter}
Assume that $H=(H^j)_{j\in J}$ satisfies \eqref{1:H0}-\eqref{1:H1} and that
\begin{equation}\label{s2:growth}
|H^j(x,r,p)-H^j(x,s,q)|\leq K(|r-s|+|p-q|) \qquad \forall r,s,p,q\in \R.
\end{equation}
Assume also that $\b_{ij}>0$ for any $i\in I_T$, $j\in Inc_i$.
Let the functions $w_1, w_2\in C^2(\G)$ satisfy \eqref{s2:CP1}.
Then $w_1\ge w_2$ on $\G$.
\end{Proposition}
\begin{Proof}
We proceed by contradiction assuming $\max_{\G} (w_2-w_1)=:\d>0$. We need the following result whose proof is postponed at the
Appendix \ref{AppendixB}.

\begin{Lemma}\label{s2:lter}
For every $\eta>0$, there exists a function $\phi_\eta\in C^2(\G)$, with $\|\phi_\eta\|_\infty\leq \eta$, such that the function $\bar w_\eta:=w_2+\phi_\eta$ satisfies
\[-\e \pd_j^2\bar w_\eta+H^j(x,\bar w_\eta,\pd_j\bar w_\eta)< -\e \pd_j^2 w_1+H^j(x,w_1,\pd_j w_1), \qquad S_\b \bar w_\eta>0.\]
\end{Lemma}

Set $\phi:=\phi_{\d/3}$ and $\bar w:= \bar w_{\d/3}$ (here, the functions $\phi_\eta$ and $\bar w_\eta$ are those introduced in Lemma \ref{s2:lter}).
We note that $\bar \d:=\max_{\G}(\bar w -w_1)>2\d/3$ and $\bar w(v_i)-w_1(v_i)\leq \d/3$ for every $i\in I_B$; therefore, for $B:=\{x\in \G:\, \bar w(x)-w_1(x)=\bar \d\}$, there holds $B\cap \G\ne \emptyset$. In fact, we claim that $B\subset \cup_{j\in J}e_j$, namely
\begin{equation}\label{s2:cl31}
v_i\notin B \qquad \forall i\in I_T.
\end{equation}
In order to prove this relation, we assume by contradiction that $v_i\in B$ for some $i\in I_T$. By Lemma~\ref{s2:lter}, we have $S^i_\b(\bar w-w_1)>0$; in particular, there exists $j\in Inc_i$ such that $\b_{ij}a_{ij}\pd_j (\bar w-w_1)>0$. This inequality contradicts the presence of a maximum at $v_i$; whence, our claim \eqref{s2:cl31} is established.

Fix $\hx\in B$. Relation~\eqref{s2:cl31} guarantees that $\hx$ belongs to some $e_j$ and that both the extremities of~$e_j$ do not belong to~$B$.
This is impossible by standard arguments; we refer the reader to \cite[Prp3.3]{koi} for a detailed proof.
\end{Proof}

\section{A priori estimates for viscous equations}\label{s3}
This section is devoted to some a priori bounds for the  the viscous equation
\begin{equation}\label{HJe2}
\left\{
\begin{array}{lr}
-\e\pd^2 w+ H^j (x,w(x),\pd  w)=0\qquad&\text{$x\in e_j$, for all $j\in J$}\\
    S^i_\b w=0    & i\in I_T\\
    w(v_i)=g_i    & i\in I_B.
\end{array}
\right.
\end{equation}
We assume that
\begin{itemize}
\item[$\bullet$]$H=(H^j)_{j\in J}$ satisfies \eqref{1:H0}-\eqref{1:H4} and either~\eqref{s2:reeegular} or~\eqref{Hr} or~\eqref{s2:growth};
\item[$\bullet$] there exist $\d>0$ and $\psi\in C^2(\G)$ such that
\begin{equation}\label{boundary}
H(x,\psi,\pd \psi)\le -\d\quad \textrm{on }\G\setminus V,\quad
S^i_\b \psi\geq 0\quad j\in I_T,\quad \psi(v_i)= g_i\quad i\in I_B;
\end{equation}

\item[$\bullet$] $\b_{ij}>0$ for every $i\in I_T$, $i\in Inc_i$.
\end{itemize}
The proof of the next two lemmas is postponed to the Appendix \ref{AppendixB}.
\begin{Lemma}\label{s3:L1}
Let $\th,\eta\in\R$, $\th>0$. Then there exists a number $M_{\th,\eta}>0$ such that
\begin{equation}\label{s3:2}
   H^j(x,r,p)>\th\quad\text{for all $p\in\R$, $|p|>M_{\th,\eta}$, $r\ge \eta$, $x\in [0,l_j]$, $j\in J$.}
\end{equation}
\end{Lemma}

\begin{Lemma}\label{s3:L2} There is a function $\phi\in C^2(\G)$ and  a vector $(\a_j)_{j\in J}$, with $\a_j\neq 0$ for all $j\in J$, for which
\[\pd_j\phi= \a_j\quad x\in e_j, j\in J,\qquad\qquad
S^i_\b\phi>0\quad i\in I_T.\]
\end{Lemma}

\begin{Theorem}\label{teo:equilipsch}
Assume that for each $\e$, there is a solution $u_\e \in C^2_{*,\b}(\G)$ of \eqref{HJe2}. Then there is $\bar \e$ sufficiently small such that
for any $0<\e<\bar \e$, the functions $u_\e$ are  uniformly bounded and equi-Lipschitz continuous on $\G$.
\end{Theorem}
\begin{Proof}\hfill\\
\noindent\textbf{ Bound on $|u_\e|$.}
For $\e$ sufficiently small, the function $\psi$ in \eqref{boundary} satisfies $\e\pd^2\psi\geq -\d$ and also
\[-\e \pd^2\psi+H(x,\psi,\pd\psi)\leq \d +H(x,\psi,\pd\psi)<0.\]
On the other hand, it fulfills $S^i_\b\psi\geq 0$ for any $i\in I_T$ and $\psi(v_i)=g_i$ for any $i\in I_B$.
By Corollary~\ref{s2:C2} (or Proposition~\ref{s2:C2bis} or Proposition~\ref{s2:C2ter}),
we get the lower bound
\begin{equation}\label{stimabasso}
\psi\leq u_\e\quad \text{on }\G,\,\text{for $\e$ sufficiently small.}
\end{equation}
To get the upper bound, we consider a function $\phi$ as in Lemma \ref{s3:L2}  and we set $\a:=\min_{j\in J}|\a_j|$. Define a function
$W\in C^2(\G)$ by $W:=M_{0,0}\phi/\a+C$, where $M_{0,0}$ as in  Lemma \ref{s3:L1}  and choose the constant $C$ in such a way that
\begin{equation*}
    W(x)>\max\{0, \max_{i\in I_B}g_i\} \quad\text{for $x\in\G$.}
\end{equation*}
By construction we have
\begin{equation}\label{s3:5}
   W^j(x)\geq 0,\quad |\pd_j W^j(x)|>M_{0,0},\quad \pd_j^2 W(x)=0  \quad\text{ for $x\in e_j$, $j\in J$.}
\end{equation}
By \eqref{s3:2} and \eqref{s3:5}, we infer
\[-\e \pd_j^2 W+H^j(x,W,\pd_j W)=H^j(x,W,\pd_j W)>0 \quad\text{ for $x\in e_j$, $j\in J$.}\]
Moreover $S^i_\b W>0$ for all $i\in I_T$ and $W(v_i)\ge g_i$  for all $i\in I_B$.
Invoking again Corollary \ref{s2:C2} (or Proposition~\ref{s2:C2bis} or Proposition~\ref{s2:C2ter}) we get the upper bound: $u_\e\le W$ on $\G$, for any $\e>0$.
We conclude that there is a constant $C_1$, independent of $\e$, such that, for $\e$ sufficiently small, there holds
\begin{equation}\label{stima}
  \max_{x\in\G} |u_\e|\le C_1.
\end{equation}

\noindent \textbf{ Bound on $|\pd^j u_\e|$}.
We split the proof in three steps devoted respectively to boundary vertices, to transition vertices and to interior of edges.

\textit{Step 1: Bound on $|\pd_j u_\e(v_i)|$, for $i\in I_B$, $j\in Inc_i$.}
Let $d_{\pd \G}:\G\to\R$ be the distance from the boundary of $\G$, i.e. $d_{\pd \G}(x):=\min\{d(x,v_i):\, i\in I_B \}$
where $d$ is the path distance on the network. For $\b>0$ set
$\G_\b:=\{x\in\G:\,d_{\pd \G}(x)\le \b\}$. We show that there are constants $K>0$, $\b>0$ and $\bar \e$ such that
\begin{equation}\label{s3:6}
\psi\le u_\e\le \psi+Kd_{\pd \G}\quad\text{on $\G_\b$, for all $0<\e<\bar \e$,}
\end{equation}
where $\psi$ is as in \eqref{boundary}.
The former inequality has been established in \eqref{stimabasso}. In order to prove the latter inequality,
%
 let $\b$ be such that $d_{\pd \G}$ does not obtain a local  maximum on the interior of $\G_\b$  and such   that there is no $i\in I_T$ for which  $v_i\in\G_\beta$. It follows that  for any  $i\in I_B$ and $j\in Inc_i$,  $|\pd_j d_{\pd\G}^j|\equiv 1$ and $|\pd^2_j d_{\pd\G}^j|\equiv 0$ on $\G_\b$.
Let
\[ \th:=\bar \e\max_{j\in J}\max_{e_j} \pd^2_j \psi^j,\quad \eta:=\min_{j\in J}\min_{e_j}\psi^j\]
and define $M_{\th,\eta}$ as in Lemma \ref{s3:L1}. Set $K:=M_{\th,\eta}+\max_{j\in J}\max_{e_j}|\pd_j\psi^j(x)|$ and
$\bar \psi:=\psi+Kd_{\pd \G}$. Hence $|\pd_j\bar\psi(x)|>M_{\th,\eta}$ for $x\in  [0,l_j]$ and by \eqref{s3:2}
\[
-\e\pd^2_j \bar\psi+H(x,\bar\psi,\pd_j \bar\psi)\ge -\th+  H(x,\bar\psi,\pd_j \bar\psi)>0,\qquad x\in\G_\b.
\]
By possible enlarging $K$, we can assume that
\[\bar\psi(x)\ge u_\e(x) \quad\text{for $x\in\pd\G_\beta\cap (\G\setminus \pd \G)$.}\]
By Corollary \ref{s2:C2} (or Proposition~\ref{s2:C2bis} or Proposition~\ref{s2:C2ter}), on each segment $e_j\cap \G_\b$ (recall that $\G_\b\cap\{v_i\}_{i\in I_T}$ is empty) we get that $\bar\psi\ge u_\e$ for any $0<\e<\bar \e$; hence relation~\eqref{s3:6} is completely proved.

By \eqref{s3:6} it follows that there exists a constant $C_2$, independent of $\e$, such that
\begin{equation}\label{s3:7}
   |\pd_ju_\e(v_i)|\le C_2\qquad \forall i\in I_B,\quad\forall 0<\e\leq \bar \e.
\end{equation}
\indent\textit{Step 2: Bound on $|\pd_j u_\e(v_i)|$, for $i\in I_T$, $j\in Inc_i$.}
We claim that  there exists a constant $C_3$ such that
\begin{equation}\label{s3:7b}
|\pd_ju_\e(v_i)|\le C_3 \quad  \forall i\in I_T,\, j\in Inc_i,\, 0<\e<\bar\e.
\end{equation}
If the claim is false, there exist $i\in I_T$, $k\in Inc_i$ and a sequence $\e_n\to 0$  such that, for $u_n:=u_{\e_n}$, we have
\begin{equation*}
    \lim_{n\to\infty}|\pd_k u_n(v_i)|=+\infty.
\end{equation*}
Let us recall: $S^i_\b u_n=\sum_{j\in Inc_i}\b_{ij} a_{ij}\pd_ju_n(v_i)=0$ for any $n\in\N$. Hence, by passing to a subsequence, there exists $j\in Inc_i$ such that
$\lim_n a_{ij}\pd_j u_n(v_i)=+\infty$.
Wlog, assume $a_{ij}=1$.
Hence, there exists a sequence $x_n\in e_j$ with $x_n\to v_i$ such that
\begin{equation}\label{s3:8}
    \lim_{n\to\infty} \pd_j u_n(x_n)=+\infty.
\end{equation}
Set $y_n:=\pi_j^{-1}(x_n)$ and fix $t_0>0$ such that $y_n+t\in [0,l_j]$ for all $t\in [0,t_0]$ and $n\in \N$.
(Note that $t_0$ is independent of $n$; indeed, as $n\to +\infty$, $y_n$ converges to $0$).
For $f_n(t):=u_n^j(y_n+t)$, relation~\eqref{s3:8} is equivalent to
\begin{equation}\label{s3:9}
    \lim_{n\to\infty} f_n'(0)=+\infty.
\end{equation}
Substituting in \eqref{HJe2} (recall: $f_n\in C^2([0,t_0])$), we get
\begin{equation}\label{s3:9b}
f''_n(t)=\e_n^{-1}H^j(y_n+t, f_n(t),f_n'(t))\quad\text{ for all $t\in [0,t_0]$, $n\in \N$}.
\end{equation}
For $C_1$ as in \eqref{stima}, set
\begin{equation}\label{s3:10}
   \th:=2C_1/t_0^2\quad\text{and}\quad \eta:=-C_1
\end{equation}
Let $M_{\th,\eta}$ be as in \eqref{s3:2}. Then
by \eqref{s3:9} there is $n\in\N$ such that $|f'_n(0)|=f'_n(0)>M_{\th,\eta}$. By \eqref{stima}, \eqref{s3:9b} and Lemma \ref{s3:L1}, we have for $\e_n<1$
\begin{equation}\label{s3:10b}
    f''_n(0)>\e^{-1}_n\th>\th.
\end{equation}
We claim that there holds
 \begin{equation}\label{s3:11}
  f_n^{''}(t)>\th\qquad  \text{for all $t\in [0,t_0]$.}
 \end{equation}
  For this purpose we set $A:=\{t\in [0, t_0]:\, f''_n(t)\ge \th\}$. By \eqref{s3:10b}
there is a connected subset $A_0$ of $A$ which contains $0$. Since $f_n\in C^2([0, t_0])$, $A_0$ is closed, hence there is a maximal $\bar t\in A_0$.
If \eqref{s3:11}  is false, then $\bar t<t_0$. Since $f'_n(0)>M_{\th,\eta}$ and $f''_n(s)\ge \th>0$ for $s\in A_0$ and therefore $f'_n$ is increasing in $A_0$, there is a neighborhood $U\subset [0,t_0]$ of $\bar t$ such that $f'_n(s)>M_{\th,\eta}$ for all $s\in U$.
Then Lemma \ref{s3:L1} and \eqref{s3:9b} imply that  $f''_n(s)>\th$ for all $s\in U$, contradicting the maximality of $\bar t$ so claim \eqref{s3:11} is proved.

Relation~\eqref{s3:11} entails the inequality
\[f_n(t)\ge \th t^2+f_n'(0)t+f_n(0)\qquad \forall t\in[0,t_0].\]
Taking into account $f'_n(0)> 0$ and \eqref{s3:10}, we estimate
\[u_n^j(y_n+t_0)=f_n(t_0)> f_n(0)+\th t_0^2\geq-C_1+\th t_0^2=C_1.\]
This relation contradicts the definition of $C_1$, hence \eqref{s3:7b} is proved.

\indent\textit{Step 3: Bound on $|\pd^j u_\e|$ on $\G$.}
By later contradiction, let us assume that $|\pd_j u_\e|$ are not uniformly bounded in $\G$, namely, there exist two sequences $\{\e_n\}_{n\in \N}$ and $\{x_n\}_{n\in \N}$, with $x_n\in\G\setminus V$, such that $|\pd_j u_{\e_n}(x_n)|\to +\infty$.
Possibly passing to a subsequence, by the compactness of $\G$, there exist $j\in J$ and $\hx\in \bar e_j$ such that $x_n\to \hx$ and $|\pd_j u_n(x_n)|\to +\infty$ for $u_n:=u_{\e_n}$.

\textit{Case (a): $\hx\in e_j$ and $\pd_j u_n(x_n)\to +\infty$.}
We shall argue as in Step~$2$; for $y_n:=\pi_j^{-1}(x_n)$, we fix $t_0>0$ such that $y_n+t\in [0,l_j]$ for all $t\in [0,t_0]$ and $n\in \N$. (Note that such a $t_0$ exists since $\hx\in e_j$).
The functions $f_n(t):=u^j_n(x_n+t)$ satisfy relations~\eqref{s3:9} and \eqref{s3:9b}.
For $\th$ and $\eta$ as in \eqref{s3:10}, we can fix $n$ sufficiently large to have $|f'_n(0)|=f'_n(0)>M_{\th,\eta}$. By \eqref{stima}, \eqref{s3:9b} and Lemma \ref{s3:L1}, we have $f''_n(0)>\th$.
We obtain relation~\eqref{s3:11} and then we conclude the proof following the same arguments as before.

\textit{Case (b): $\hx\in e_j$ and $\pd_j u_n(x_n)\to -\infty$.}
We shall use arguments analogous to those of previous case. Fix $t_0>0$ such that $y_n-t\in [0,l_j]$ for all $t\in [0,t_0]$ and $n\in \N$. (Note that such a $t_0$ exists since $\hx\in e_j$).
The functions $f_n(t):=u_n^j(y_n-t)$  satisfy relation~\eqref{s3:9} and
\begin{equation}\label{s3:9ter}
f''_n(t)=\e_n^{-1}H^j(y_n-t, f_n(t),-f_n'(t))\quad\text{ for all $t\in [0,t_0]$, $n\in \N$}.
\end{equation}
Fix $\th$ and $\eta$ as in \eqref{s3:10}; fix $n$ sufficiently large to have $-f'_n(0)<-M_{\th,\eta}$. By \eqref{stima}, \eqref{s3:9ter} and Lemma \ref{s3:L1}, we have $f''_n(0)>\th$.
We obtain relation~\eqref{s3:11} and then we conclude the proof following the same arguments as before.
%
%

\textit{Case (c): $\hx=v_i\in V$ and $\pd_j u_n(x_n)\to -\infty$.}
Wlog, we assume $a_{ij}=1$ (recall that $e_j$ is the edge containing all the $x_n$).
Fix $n$ sufficiently large to have
\begin{equation}\label{s3:15}
\pd_j u_n(x_n)<-\max\{C_2,C_3,\bar C\}
\end{equation}
where $C_2$ and $C_3$ are respectively the constant introduced in \eqref{s3:7} and in \eqref{s3:7b} while $\bar C$ is such that
\begin{equation}\label{s3:16}
H(x,-C_1,p)> 0\qquad \forall x\in\G, |p|>\bar C
\end{equation}
(assumption \eqref{1:H4} ensures the existence of the constant $\bar C$).
For each $n\in \N$, let $t_n\in(0,l_j)$ be such that $y_n-t\in [0,l_j]$ for all $t\in [0,t_n]$. Observe that in this case $t_n$ depends on $n$ and that $\pi_j(y_n-t_n)=v_i$.
By assumption~\eqref{1:H1}, for every $n\in\N$, the function $f_n(t):=u_n^j(y_n-t)$ satisfies relation~\eqref{s3:9} and also
\begin{equation}\label{s3:9quat}
f''_n(t)=\e_n^{-1}H^j(x_n-t, f_n(t),-f_n'(t))\geq
\e_n^{-1}H^j(x_n-t, -C_1,-f_n'(t))
\end{equation}
for every $ t\in [0,t_n]$. Taking into account relations \eqref{stima}, \eqref{s3:15}, \eqref{s3:16} and \eqref{s3:9quat}, we infer: $f''_n(0)>0$. In fact, let us prove
 \begin{equation}\label{s3:11quat}
f_n^{''}(t)>0 \qquad  \forall t\in [0,t_n].
\end{equation}
In order to prove this inequality, we introduce the set $A:=\{t\in [0,t_n]:\, f''_n(t)\geq 0\}$ and the set $A_0$ as its connected component containing $t=0$. Let $\bar t$ be the maximal point of $A_0$; for later contradiction, assume that $\bar t<t_n$. We observe the function $f'_n$ is increasing in $(0,\bar t)$ and, by \eqref{s3:15}, $f'_n(0)>\max\{C_2,C_3,\bar C\}$. Hence, it follows that: $f'_n(\bar t)>\max\{C_2,C_3,\bar C\}$ and, by \eqref{s3:9quat}, $f''_n(\bar t)>0$. A contradiction to the maximality of $\bar t$ is obtained so inequality \eqref{s3:11quat} is completely proved.

Relations \eqref{s3:15} and \eqref{s3:11quat} entail
\[\pd_j u_n(v_i)=-f'_n(t_n)<-f'(0)=\pd_j u_n(x_n)<-\max\{C_2,C_3,\bar C\}\]
which contradicts the definition either of $C_2$ or of $C_3$.

\textit{Case (d): $\hx=v_i\in V$ and $\pd_j u_n(x_n)\to +\infty$.}
In this case, it suffices to follow the same arguments of Step~$2$.
\end{Proof}
\begin{Remark}
This theorem applies to problem \eqref{HJe}. In fact,  a priori estimates for this problem could be obtained by  \cite[Thm2, App1]{li}. However, for the sake of completeness, a direct proof has been given.
\end{Remark}
\section{The vanishing viscosity limit}\label{s4}
In this section we prove the vanishing viscosity result, i.e. the convergence of the solution of \eqref{HJe2} to the one of \eqref{HJ}. We observe that assumptions \eqref{1:H5}-\eqref{1:H6} are not necessary for \eqref{HJe2} but they play an crucial role for the uniqueness of \eqref{HJe2}. Moreover the specific form of the Hamiltonian in \eqref{HJe} is only used to prove the existence of a solution, while a priori estimates in section \ref{s3} and the convergence of the vanishing viscosity limit in this section hold for the more general class of Hamiltonians.
\begin{Theorem}\label{s4:T1}
 Assume that $H=(H^j)_{j\in J}$ satisfies \eqref{1:H0}-\eqref{1:H6}.
Let $u_n:=u_{\e_n}\in C^2_{*,\b}(\G)$ be a sequence of solutions of \eqref{HJe2} such that $u_n$ and $\pd  u_n$ are uniformly bounded on $\G$.
If  $u_n$ converges uniformly to a function $u\in C(\G)$,
then $u$ is a solution of \eqref{HJ}.
\end{Theorem}

For the proof  we need two lemmas: the former is an immediate consequence of \eqref{1:H0}--\eqref{1:H6} while the proof of the latter is postponed to Appendix \ref{AppendixB}.
\begin{Lemma}\label{s4:L1}
Under the hypotheses of Theorem~\ref{s4:T1}, for $i\in I_T$, define a function $h_i:\R\to\R$ by $h_i(p):=H^j(v_i,0,p)$, $j\in Inc_i$ (by \eqref{1:H5} the definition is independent of $j$). Then,
 $h_i(0)=\min h_i$, $h_i$ is symmetric and nondecreasing on $(0,+\infty)$. In particular, either it is strictly positive or there there is a unique number $a\geq 0$ such that $h(a)=h(-a)=0$.
\end{Lemma}
\begin{Lemma}\label{s4:L2}
Assume the hypotheses of Theorem~\ref{s4:T1}.
Let $i\in I_T$, $j\in Inc_i$ and $\xi>0$. Furthermore let $x_m\in e_j$, $m\in \N$, such that $\lim_m x_m=v_i$. Then there is a number
$m_\xi\in\N$ such that for all $m>m_\xi$
 \begin{equation}\label{prp:noconvex}
    H^j\left(v_i, u(v_i), \frac{u(x_m)-u(v_i)}{d(x_m,v_i)}\right)\le \xi.
\end{equation}

\end{Lemma}

\begin{Proofc}{Proof of Theorem \ref{s4:T1}}\\
{\it Step 1: $u$ is a subsolution of \eqref{HJ}.}
For $x\in e_j$ (for some $j\in J$), the proof is standard and we skip it (see \cite[Thm2.3]{b}).
Assume that $x=v_i$, for some $i\in I_T$.
Let $j,k\in Inc_i$, $j\neq k$ and let $\phi$ be  a $(j,k)$-test function such that  of $u-\phi$ has a local maximum  at $x$.
We shall assume $u(0)=0$; the general case can be dealt with by similar arguments and we shall omit it.
Then
\[h(\pd_j\phi(x))=H^j(x, u(x),\pd_j\phi(x))\]
where $h=h_i$ as in Lemma \ref{s4:L1}.
We claim that $h$ is not strictly positive; actually, by contradiction, let us assume $h>0$. In particular, we have $h(0)>0$ and, by the continuity of $H^j$, we infer $H(x,u(x),0)>0$ in some $B_\eta(v_i)$. By Lemma~\ref{s4:L1}, we get $H(x,u(x),\p_j\phi(x))>0$ for every test function at some points in $B_\eta(v_i)$. This inequality contradicts that $u$ is a subsolution in $e_j$.
By Lemma~\ref{s4:L1}, there exists a unique number $a$ such that $h(p)>0$ for $|p|>a$.

Suppose by contradiction that $h(\pd_j\phi(x))>0$. Since $\phi$ is $(j,k)$-differentiable
at $x$ and therefore $a_{ij}\pd_j \phi(v_i)+a_{ik}\pd_k \phi(v_i)=0$, for one of the indices $j,k$, say for $j$, there is a number $\d_0>0$ such that
\begin{equation}\label{s4:2}
    a_{ij}\pd_j \phi(v_i  )=-(a+\d_0)
\end{equation}
where $a>0$ is defined as in Lemma \ref{s4:L1}. Let $x_m$ be a sequence with $x_m\in e_j$ with $\lim_{m\to \infty} x_m=x$. As $u-\phi$  attains a local maximum at $x$, by \eqref{s4:2} we get 
\[p_m:=\frac{u(x_m)-u(x)}{d(x_m,x)}\leq \frac{\phi(x_m)-\phi(x)}{d(x_m,x)}<-(a+\frac{\d_0}{ 2})\]
for $m$ sufficiently large. By the properties of $h$ it follows that there exists $\d_1>0$ such that
\[\d_1<h(p_m)=H^j(x,u(x),p_m)\]
for $mù$ sufficiently large, a contradiction to Lemma \ref{s4:L2}. Hence
\[h(\pd_j\phi(x))=H^j(x, u(x),\pd_j\phi(x))\le 0.\]
{\it Step 2: $u$ is a supersolution of \eqref{HJ}.}
For $x\in e_j$ (for some $j\in J$), the proof is standard and we skip it (see \cite[Thm2.3]{b}). Assume that $x=v_i$, for some $i\in I_T$.
The proof is based on the following lemma (the proof is in Appendix \ref{AppendixB}).
\begin{Lemma}\label{s4:L3}
Assume the hypotheses of Theorem~\ref{s4:T1}.
Let $i\in I_T$ and assume that, for $j\in Inc_i$, there holds $a_{ij}\pd_j u_n(v_i)\le 0$ for infinitely many $n\in \N$.
Furthermore assume that there is a function $\phi\in C^2(\G)$ such that $u-\phi$ has a local minimum at $v_i$. Then
$H^j(v_i,u(v_i),\pd^j\phi(v_i))\ge 0$.
\end{Lemma}
Since   $u_n$ satisfies \eqref{Kir} at $x=v_i$, there is an index $j\in Inc_i$ such that
\begin{equation}\label{s4:7}
    a_{ij}\pd_j u_n(v_i)\le 0
\end{equation}
for infinite many $n\in \N$.
 We show that $j$ is a $k$-feasible index for each $k\in Inc_i\setminus\{j\}$.
 We assume wlog that $a_{ij}=1$  and we fix a $(j,k)$-test function  $\phi$ such that $u -\phi$   has a strict minimum point at $0=\pi_j^{-1}(v_i)$ relatively to $\bar e_j\cup \bar e_k$. Let $\phi_m\in C^2([0,l_j])$ ($m\in\N$), be such that $\phi_m$ converges to $\phi$ with respect the topology of $C^1([0,l_j])$. Let $z_m\in \bar e_j\cup \bar e_k$ be such that $u-\phi_m$ attains a local minimum with respect to $\bar e_j\cup \bar e_k$. Then, by standard arguments, the point $z_m$ converges to $x$ and either by the case $x\in e_j$ if $z_m\in e_j$ or by Lemma  \ref{s4:L3} if $z_m=x$, we conclude that
  \[H^j(z_m,u(z_m),\pd^j\phi_m(z_m))\ge 0.\]
 Since $\lim_{m\to\infty}\pd_j\phi_m(z_m)=\pd_j\phi(x)$, we obtain
\[H^j(v_i,u(v_i),\pd^j\phi(v_i))
=\lim_{m\to\infty} H^j(z_m,u(z_m),\pd^j\phi_m(z_m))\ge 0.\]
Hence $j$ is $i$-feasible for $k$ and by symmetry $k$ is $i$-feasible for $j$ at $x$.
\end{Proofc}
\subsection{Example: the eikonal equation}\label{esempio}
We consider the eikonal equation on the network $\G$ with null boundary condition
\begin{equation}\label{E1}
|\pd u|=f(x) \qquad \textrm{on }\G, \qquad\qquad u(v_i)=0\qquad \forall i\in I_B
\end{equation}
where $f$ is a Lipschitz continuous function with $f\geq \a>0$.

\noindent\texttt{Fact 1.} {\it There exists a unique viscosity solution~$u$ to~\eqref{E1}.}

For the proof, we refer the reader to \cite{sc} (see also \cite{csm} for the generalization to LEP spaces); in fact, $u$ can be written as a {\it weighted} distance from $\pd \G$.

We observe that a function $u$ solves \eqref{E1} if, and only if, it solves
\begin{equation}\label{E2}
|\pd u|^2=f^2(x) \qquad \textrm{on }\G, \qquad\qquad u(v_i)=0\qquad \forall i\in I_B.
\end{equation}
For any collection $\b=(\b_{ij})$ ($i\in I_T$, $j\in Inc_i$) with $\b_{ij}>0$, we introduce the viscous approximation to \eqref{E2}:
\begin{equation}\label{E3}
-\e \pd^2u+|\pd u|^2=f^2(x) \quad \textrm{on }\G
, \quad S^i_\b u=0\quad \forall i\in I_T
, \quad u(v_i)=0\quad \forall i\in I_B.
\end{equation}
\noindent\texttt{Fact 2.} {\it By Theorem~\ref{Teo:logar}, there exists a unique classical solution~$u_\e$ to~\eqref{E3}.}

\noindent\texttt{Fact 3.} {\it By Theorem~\ref{teo:equilipsch}, the functions $u_\e$ are equibounded and equilipschitz continuous.}

\noindent\texttt{Fact 4.} {\it The sequence $\{u_\e\}$ uniformly converges to $u$.}

Actually, by facts 3, Ascoli's Theorem ensures that there exists a subsequence $\{u_{\e_n}\}$ uniformly convergent to some function~$v$. By Theorem~\ref{s4:T1}, $v$ is a solution to \eqref{E2}.  By the uniqueness of the solution to \eqref{E2}, we deduce that the whole sequence $\{u_\e\}$ converges to its unique solution $u$.



\appendix
\section{Appendix}\label{AppendixB}
\begin{Proofc}{Proof of Lemma \ref{s2:lter}}
Fix two functions $w_1,w_2\in C^2(\G)$ such that relations \eqref{s2:CP1} hold. By the regularity of $w_1$, we can introduce $\tilde H^j(x,r,p):=H^j(x,r,p) +\e \pd_j^2 w_1-H(x,w_1,\pd_j w_1)$. For $\bar w_\eta:=w_2+\phi_\eta$, assumptions \eqref{s2:CP1} and \eqref{s2:growth} entail
\[-\e\pd_j^2 \bar w_\eta+\tilde H ^j(x,\bar w_\eta,\pd_j \bar w_\eta)\leq
-\e \pd_j^2\phi_\eta+K(\|\phi^j_\eta\|_\infty+\|\pd_j\phi^j_\eta\|_\infty ).\]
Therefore, it is enough to prove that, for every $\eta>0$ there exists $\phi_\eta$ such that
\begin{equation}\label{claim:L3.1}
 \|\phi_\eta\|\leq \eta,\quad S_\b \phi_\eta>0,\quad -\e \pd_j^2\phi_\eta+K(\|\phi^j_\eta\|_\infty+\|\pd_j\phi^j_\eta\|_\infty )<0.
\end{equation}
Let $\d:I\times I\to\N$ be the metric given by the smallest number $\d(i,j)$ of the edges   a path connecting $v_i$ and $v_j$ can consist of. It induces a partition $I_l:=\{i\in I:\, \d(i,I_B)=l\}$. Observe that $I_0=I_B$ and set $m:=\max\{l\in\N:\, I_l\neq \emptyset\}$.

For simplicity, we address only the case $m=1$ with $l_j=l$ for $j\in J$; the general case can be dealt with in a similar manner and we shall omit it. In this case, each vertex belongs either to $\G_0:=\{v_i: i\in I_0\}$ or to $\G_1:=\{v_i: i\in I_1\}$; furthermore, by Remark~\ref{rmk-A}, each edge connects either two vertices in $\G_1$ or a vertex in $\G_0$ and one in $\G_1$ (namely, it do not connect two vertices in $\G_0$).

Let us enumerate the elements in $\G_1$ as $\{v_{i_1},\dots, v_{i_n}\}$. Wlog, we assume that: when $e_j$ connects $v_{i_s},v_{i_t}\in \G_1$, with $1\leq s<t\leq n$, its parametrization is $\pi_j(0)=v_{i_s}$, $\pi_j(l_j)=v_{i_t}$ while, for $e_j$ connecting $v_{i_s}\in \G_1$ and $v_k\in \G_0$, its parametrization is $\pi_j(0)=v_{i_s}$, $\pi_j(l_j)=v_k$.
Let us now define a function $\phi\in C(\G)$ in the following manner: on the vertices, we set
\begin{equation*}
\phi(v_{i_s}):=e^{2K(s-1)l\e^{-1}}\quad \forall v_{i_s}\in \G_1,\qquad
\phi(v_k):=e^{2K(n+1)^2 \b_0 l\e^{-1}}\quad \forall v_k\in \G_0
\end{equation*}
with $\b_0:= \max \b_{ij}/\min \b_{ij}$; moreover, on the edge $e_j$, we set
\begin{align*}
&\phi^j(x):=e^{2K(s-1)l\e^{-1}}e^{2K(t-s)\e^{-1} x}&&\qquad \textrm{if $e_j$ connects $v_{i_s}$ and $v_{i_t}$, $s<t$}\\
&\phi^j(x):=e^{2K(s-1)l\e^{-1}}e^{2K[(n+1)^2\b_0-s+1]\e^{-1}x}&&\qquad \textrm{if $e_j$ connects $v_{i_s}\in \G_1$ and $v_k\in \G_0$}.
\end{align*}
One can easily check that, on each edge $e_j$, last relation of~\eqref{claim:L3.1} is satisfied.
On the other hand, for $J_1:=\{j\in Inc_{i_s}:e_j\textrm{ connects }v_{i_s} \textrm{ with some } v_{i_t}\in \G_1\}$ and $J_2:=\{j\in Inc_{i_s}:e_j\textrm{ connects }v_{i_s} \textrm{ with some } v_k\in \G_0\}$, we have
\[S^{i_s}_\b \phi=\sum_{j\in J_1}\b_{i_s j} a_{i_sj} \pd_j\phi(v_{i_s})+\sum_{j\in J_2}\b_{i_sj} a_{i_sj} \pd_j\phi(v_{i_s})\equiv S_1+S_2.\]
Since $\# J_2\geq 1$ and $a_{i_sj}=1$ for $j\in J_2$, we infer
\[S_2\geq 2K(\min \b_{ij})e^{2K(s-1)l}[(n+1)^2\b_0-s+1]\e^{-1}\geq
2K(\max \b_{ij})e^{2K(s-1)l}[(n+1)^2-n]\e^{-1}.\]
On the other hand, since $\#J_1\leq n-1$, we get
\[S_1\geq -\sum_{t=1}^n \b_{i_sj}e^{2K(s-1)l}(t-s)\e^{-1}\geq
-K(\max \b_{ij})e^{2K(s-1)l}(n+1)^2\e^{-1}.\]
Owing to the last three relations, we have $S^{i_s}_\b\phi>0$ for $s=1,\dots,n$.

Finally, we observe that relations in~\eqref{claim:L3.1} are linear; whence, the function $\phi_\eta:=\eta \frac{\phi}{\|\phi\|}$ is a desired function.
\end{Proofc}
\begin{Proofc}{Proof of Lemma  \ref{s3:L1}}
Fix $\theta$ and $\eta$ as in the statement. By~\eqref{1:H1}, we have: $H^j(x,r,p)\geq H^j(x,\eta,p)$  for every $x\in e_j$, $r\geq \eta$, $p\in\R$, $j\in J$.
By~\eqref{1:H4}, there exists $M_{\th,\eta} >0$ such that:
$H^j(x,\eta,p)>\theta$ for every $x\in e_j$, $r\geq \eta$, $|p|>M_{\th,\eta}$, $j\in J$.
Substituting the previous inequality in the last one, we accomplish the proof.
\end{Proofc}
\begin{Proofc}{Proof of Lemma  \ref{s3:L2}}
Define the set
\begin{equation*}
    M:=\{\xi\in\R^I:\,\text{$\xi_i\neq\xi_j$ for all $i,j\in I$ with $v_i\, \mathrm{adj}\,v_j$}\}
\end{equation*}
and observe that there is an injective map $\Phi:M\to D$ with
\[  D:=\{\phi\in C^2(\G):\, \text{there exists $(\a_j)_{j\in J}$ s.t. $\a_j\neq 0$ and $\pd_j\phi\equiv \a_j$ on $e_j$, $j\in J$}\}\]
such that $\Phi[\xi](v_i)=\xi_i$, $i\in I$. It suffices to show that there is a $\xi\in M$ such that $S^i_\b(\Phi[\xi])>0$ for all $i\in I_T$.
To this end, we define $I_l$ and $m$  as in the proof of Lemma~\ref{s2:lter} and, for $i\in I_T$, we introduce the map
$T_i:=S^i_\b\circ \Phi$
which is: (a) continuous, unbounded and strictly decreasing in the component $\xi_i$, (b) continuous, unbounded and strictly increasing in each component $\xi_j$, $j\in A_i:=\{j\in I:\, v_j \mathrm{adj}\,v_i\}$, (c) independent of $\xi_j$ for any $j\in I\setminus ( \{i\}\cup A_i)$.

Let us now construct $\xi\in M$ such that $T_i(\xi)>0$ for all $i\in I_T$.
We first choose $\xi\in M$. Let $i\in I_m$, by property $(b)$ and by $I_{m-1}\cap A_i\neq \emptyset$, we may increase the value of $\xi_j$, $j\in I_{m-1}$, such that we obtain $T_i(\xi)>0$ for all $i\in I_m$ and such that  $\xi$ remains in $M$. Analogously, we can increase $\xi_j$, $j\in I_{m-2}$, such that $T_i(\xi)>0$ for all $i\in I_{m-1}$ and such that $\xi_j$, $j\in J$, remain pairwise different. For $k=3,\dots, m$ we continue this procedure by sufficiently increasing $\xi_j$, $j\in I_{m-k}$, in order to ensure that $T_i(\xi)>0$ for all $i\in I_{m-k+1}$, ending up with a choice for $\xi\in M$ such that $T_i(\xi)>0$ for all $i\in \cup_{l=1}^m I_l=I_T$.
\end{Proofc}
\begin{Proofc}{Proof of Lemma  \ref{s4:L2}}
Let us recall that our hypotheses entail: $\|u_n\|_\infty \leq C_1$, $\|\p u_n\|_\infty\leq C_2$, $\|u_n-u\|_\infty\to 0$, $\e_n\to 0$ and $u$ is Lipschitz continuous with a Lipschitz constant not greater than $C_2$.
For the sake of clarity, we split the proof in several steps.

{\it Step 1.}
For $\e>0$, introduce $N^x_\theta$ as
\[H^j(x, u(x)-\e,p)>\theta \qquad \forall p,\, N^x_\theta<|p|<C_2.\]
We observe that $N^{v_i}_\theta$ is non decreasing in $\theta$ by \eqref{1:H3}  and, by the continuity of $u$, there holds
\[|H^j(x,u_n(x)-\e,p)-H^j(v_i,u(v_i)-\e,p)|\leq \omega(|x-v_i|)\qquad\forall x\in \G, \, |p|<C_2\]
where $\omega$ is the modulus of continuity of $H$ on $\G\times [-2C_1,2C_1]\times [-C_2,C_2]$.
Consider $\eta>0$ such that $\omega(|x-v_i|)<\e$ and $|u(v_i)-u_n(x)|<\e$ for $x\in[0,\eta)$.
Fix $\bar x\in(0,\eta/2)$ and $\bar \eta<\bar x$; our claim is to prove that, there exist $N\in\N$ such that
\begin{equation}\label{cl:bdd_n}
|H^j(v_i,u(v_i),\p_j u_n(x))|<2\e\qquad \forall x\in(\bar x-\bar\eta,\bar x+\bar\eta),\, n>N.
\end{equation}
In order to prove this relation, we proceed by contradiction assuming that, for some $x\in (\bar x-\bar \eta,\bar x+\bar \eta)$  there holds
$H^j(v_i,u(v_i),\p_j u_n(x))\geq  2\e$ for every $n>N_0$.
By assumption~\eqref{1:H1} and the equation in~\eqref{HJe2}, we deduce
\[\e_n \p^2_j u_n(x)=H^j(x,u_n(x),\p_ju_n(x))\geq
H^j(v_i,u(v_i)-\e,\p_ju_n(x))-\e>\e \qquad \forall n>N.
\]
Therefore, we have: $u_n(x)\geq u(\bar x)-\e$, $\p_j u_n(x)\geq N^{v_i}_\e>0$ and $\e_n \p^2_j u_n(x)>\e$.
We claim that, for $N$ sufficiently large (it suffices to have $N>N_0$ and $\e_n< \e C^{-1}_2(\bar x+\bar\eta -x)$ for $n>N$), these inequalities still hold in $[x,\bar x+\bar \eta]$, namely
\begin{equation}\label{cl_cre}
u_n(y)\geq u(\bar x)-\e,\quad \p_j u_n(y)\geq N^{v_i}_\e,\quad \e_n\p^2_j u_n(y)\geq\e\qquad \forall y\in[x,\bar x+\bar\eta].
\end{equation}
Indeed, let $A$ and $\bar t$ be respectively the connect set containing $x$ where they hold and its maximum point. If $\bar t<\bar x+\bar\eta$, since $u_n$ and $\p u_n$ are both strictly increasing on $[x,\bar t]$, we have
$u_n(\bar t)> u(\bar x)-\e$, $\p_j u_n(\bar t)> N^{v_i}_\e$; by~\eqref{HJe2}, we get
\begin{equation*}
\e_n \p^2_j u_n(\bar t)=H^j(\bar t,u_n(\bar t),\p_ju_n(\bar t))
\geq H^j(v_i,u(v_i)-\e,\p_ju_n(\bar t))-\e>\e.
\end{equation*}
Hence by continuity there is a neighborhood of $\bar t$ contained in $A$; this fact contradicts the definition of $\bar t$. Claim~\eqref{cl_cre} is completely proved.

Relations~\eqref{cl_cre} and our choice of $\e_n$ ensure the following relation
\[\p_j u_n(\bar x+\bar\eta)\geq \p_j u_n(x)+\e \e_n^{-1}(\bar x+\eta-x)>C_2\]
which contradicts our bound on $\p u_n$. Hence, we get:
$H^j(x,u_n(x),\p_j u_n(x))\leq 2\e$.

In order to prove the other inequality of \eqref{cl:bdd_n} we proceed in a similar manner. We assume by contradiction: $H^j(x,u_n(x),\p_j u_n(x))\leq -2\e$. for every $n>N_0$ (for some $x\in(\bar x-\bar\eta,\bar x+\bar\eta)$). We choose $N$ such that: $N>N_0$ and $\e_n< \e C^{-1}_2(x+\bar x-\bar\eta)$ for $n>N$. Arguing as before, we infer:
\[u_n(y)\geq u(\bar x)-\e,\quad \p_j u_n(y)\leq -N^{v_i}_{-\e},\quad \e_n\p^2_j u_n(y)\geq\e\qquad \forall y\in[\bar x-\bar \eta,\bar x].\]
These relations and our choice of $\e_n$ ensure: $\p_j u_n(\bar x-\bar\eta)\leq \p_j u_n(x)-\e \e_n^{-1}(x+\bar x-\bar\eta)<-C_2$, a contradiction of our bound on $\p u_n$. Hence the bound \eqref{cl:bdd_n} is completely proved.

{\it Step 2.}
Assume wlog $a_{ij}=1$. The aim is to prove that, for each $\xi>0$, there exist $ \eta>0$ such that
\begin{equation}\label{cl:bd_vert}
H^j\left(v_i, u(v_i),\frac{u^j(0)-u^j(y)}{y}\right)\leq \xi \qquad \forall y\in(0,\eta).
\end{equation}
In order to prove this relation, for each $\e>0$, consider $\eta$ as before.
Fix $y\in (0,\eta/2]$ and $x:=y\e C_2^{-1}$. By the Lipschitz continuity of $u$, our choice of $x$ and the uniform convergence, for $n$ sufficiently large, we infer
\begin{eqnarray*}
\frac{|u^j(0)-u^j(y)|}{y} &\leq&
\frac{|u^j(x)-u^j(y)|}{y}+C_2\frac{x}{y}\\
&\leq&  \frac{|u^j(x)-u^j(y)|}{y-x}+
\frac{|u^j(x)-u^j(y)|}{y-x}\left(\frac {y-x}y-1\right)+C_2\frac xy
\\ &\leq& \frac{|u^j_n(x)-u^j_n(y)|}{y-x}+2\frac{\|u_n-u\|_\infty}{y-x}+2\e.
\end{eqnarray*}
By mean value theorem, we deduce for any $n\in \N$
\begin{eqnarray*}
&&|H^j(v_i,u^j(v_i),\frac{u^j(0)-u^j(y)}{y})|\leq
|H^j(v_i,u^j(v_i),\frac{u^j_n(y)-u^j_n(x)}{y-x})| +\omega(2\frac{\|u_n-u\|_\infty}{|y-x|}+2\e)\\
&&\qquad\qquad\leq|H^j(v_i,u^j(v_i),\p_ju_n(x'_n))| +\omega(2\frac{\|u_n-u\|_\infty}{|y-x|}+2\e)
\end{eqnarray*}
for some $x'_n\in(x,y)$.
 Letting $n\to+\infty$, by step 1, we infer
\[|H^j(v_i,u^j(v_i),\frac{u^j(0)-u^j(y)}{y})|\leq \e +\omega(\e)\]
In conclusion, it suffices to choose $\e$ such that $\e +\omega(\e)<\xi$.
\end{Proofc}

\begin{Proofc}{Proof of Lemma  \ref{s4:L3}}
Wlog, we assume that $u(v_i)=\phi(v_i)=0$  and $a_{ij}=1$. By the assumptions, we can choose a subsequence of $(u_n)_{n\in\N}$ (still denoted by $(u_n)_{n\in\N}$) such that $\pd_ju_n(0)\le 0$ for all $n\in\N$.
Our aim is to prove that
\[h(\p\phi(x)):=h_i(\p\phi(x))\geq 0\]
where $h_i$ is the function introduced in Lemma~\ref{s4:L1}. For $h(p)\geq 0$ for every $p$, there is nothing to prove. By Lemma~\ref{s4:L1}, let us assume that there exists $a>0$ such that $h(p)<0$ on $(-a,a)$. We want to show that $\pd_j\phi(v_i)\le -a$. To this end we assume the contrary, i.e.
 there is $\d\in (0,2a)$
such that
\begin{equation}\label{a:41}
   \pd_j\phi(v_i)= -a+\d
\end{equation}
and we set $H^j(v_i,u(v_i),\pd_j\phi(v_i))=-\a<0$.
We claim that for   $n\in\N$ sufficiently large, there is $r_n>0$ such that
\begin{equation}\label{claimB2}
u_n^j(x)<u_n^j(0),\qquad \pd^ju_n (x)<0\quad \textrm{for }x\in (0,r_n].
\end{equation}
This is clear if $\pd_j u_n(v_i)<0$. Assume $\p_j u_n(v_i)=0$.
In order to prove \eqref{claimB2}, it is enough to prove that, for $n$ sufficiently large, there exists $r_n>0$ such that
\[\pd_j^2 u_n(x)<-\a/2\qquad\forall x\in(0,r_n].\]
To this end, we argue by contradiction and we assume that there exists a sequence $x_m\in e_j$, with $x_m\to v_i$ as $m\to+\infty$, such that $\pd_j^2 u_n(x_m)>-\a/2$. The continuity of $\pd_j u_n$ ensures: $\pd_j u_n(x_m)\to 0$ as $m\to +\infty$. Moreover, we have
\begin{eqnarray*}
\e_n\pd_j^2 u_n(x_m)&=&H^j(x_m, u_n(x_m),\pd_j u_n(x_m))\\&=&
H^j(v_i, u(v_i),0)+\omega(|x_m-v_i|+|u_n(x_m)-u(v_i)|+|\pd_j u_n(x_m)|)
\end{eqnarray*}
where $\omega$ is the modulus of continuity of $H$ in $\G\times[-C,C]\times[-C,C]$ and $C$ is a constant such that $\|u_n\|_\infty,\|\pd u_n\|_\infty \leq C$ (its existence is ensured by the hypotheses of Theorem~\ref{s4:T1}).
Owing to its monotonicity in $|p|$, $H^j$ fulfills
\[H^j(v_i,u(v_i),0)\le H^j(v_i,u(v_i),\pd_j\phi(v_i))\le -\a.\]
Taking into account the last two relations, we infer
\[\e_n\pd_j^2 u_n(x_m)\le -\a+\omega(|x_m-v_i|+|u_n(x_m)-u(v_i)|+|\pd_j u_n(x_m)|)\]
which gives the desired contradiction for $n$ sufficiently large; hence, \eqref{claimB2} is proved.

Let us now show that there exists $r>0$ such that, for $n$ sufficiently large, $u_n^j$ cannot obtain a local minimum in $(0,r]$. In fact, if $u_n^j$ has a minimum at $x$, we get by \eqref{1:H0} and the uniform
bound on $\pd_j u_n$
\begin{eqnarray*}
0&\le& \e_n\pd_j^2 u_n(x)=H^j(x,u_n^j(x),0) \le H^j(v_i,u_n^j(v_i),0)+\omega(|v_i-x|)\\
&\le& -\a+\o(1/n)+\omega(|v_i-x|)<0
\end{eqnarray*}
for $n$ sufficiently large and $|v_i-x|$ small,  hence a contradiction.
Therefore
$u_n^j(x)\le u_n^j(0)$ for $x\in [0,r]$. It follows that
\begin{equation}\label{a:42}
  u^j(y)=\lim_{n\to\infty} u_n^j(y)\le \lim_{n\to\infty} u_n^j(0)=u(v_i)=0\qquad \forall y\in [0,r]
\end{equation}
namely, $u^j$ attains in $0=\pi_j^{-1}(v_i)$ its maximum with respect to $[0,r]$.
Since $u-\phi$ attains a local minimum at $v_i$, \eqref{a:42} implies that we may restrict to consider the case $\d\le a$ in \eqref{a:41}.

By the continuity of $H^j$, \eqref{1:H1} and Lemma \ref{s3:L1}, it follows that there are $\eta,\g>0$ with $\eta<\min\{\d,r\}$ such that
\begin{equation}\label{a:44}
    H^j(x,z,p)\le -\g\quad\text{ for all $p\in [-\b,0]$, $z\in (-\infty,\eta]$  and $x\in [0,\eta]$}
\end{equation}
where $\b:=a-\d+\eta$. Choose $n_0$ such that $\e_{n_0}\b/\g<r$ and $u_n^j(0)<\eta$ for all $n\ge n_0$. For $n\ge n_0$   set $v_n(x):=\pd_ju_n(x)$ for $x\in (0,r)$.
By \eqref{HJe2}, \eqref{claimB2}, $u_n^j(x)\le u_n^j(0)$ for $x\in [0,r]$  and \eqref{a:44}, we get
\begin{equation}\label{a:45}
   \pd_jv_n(x)=H^j(x,u_n(x),v_n(x))/\e_n\le -\g/\e_n
\end{equation}
for all $x\in [0,\eta)$ and $-\b\le v_n(x)\le 0$. In particular, since we have $-a+\d\leq v_n(0)\le 0$, we derive from \eqref{a:45}
that there is $x_n$ with
\begin{equation}\label{a:46}
    0\le x_n\le \e_{n}\b/\g \le\e_{n_0}\b/\g<r
\end{equation}
such that $v_n(x_n)=-\b$.
We furthermore claim that
\begin{equation}\label{a:48}
    v_n(x)\le -\b \quad\text{ for all $x_n<x\le \eta$.}
\end{equation}
Actually, if the claim were not true, there would be $x_0$ with $x_n<x_0<\eta$ such that $v_n(x_0)=-\b$ and $\pd_j v_n(x_0)\ge 0$. This contradicts \eqref{a:45}.

Now, \eqref{a:48} and $u_n^j(x)\le u_n^j(0)$ for $x\in e_j$, $n\in \N$ imply
\begin{align*}
    u_n^j(y)=u_n^j(x_n)+\int_{x_n}^y v_n(s)ds\le u_n^j(x_n)-\b(y-x_n)\le u_n^j(0)-\b(y-x_n)
\end{align*}
for all $y$ with $x_n\le y\le\eta$. Using \eqref{a:46} we conclude
\[u^j(y)=\lim_{n\to\infty}u_n^j(y)\le -y\b=y(-a+\d-\eta)\qquad\forall y\in[0,\eta].
\]
As $u^j-\phi^j$ has  a local minimum at $0=\pi_j^{-1}(v_i)$, it follows that there is $\rho>0$ such that $\phi^j(y)\le y(-a+\d-\eta)$ for all $0\le y\le \rho$, a contradiction to \eqref{a:41}.
 \end{Proofc}

\noindent{\bf Acknowledgment.}
The first and the second authors have been partially supported by the Indam project ``{\it Fenomeni di propagazione sui grafi ed in mezzi eterogenei}''.


\end{document}